\documentclass[11pt]{amsart}

\usepackage[
  margin=1in
]{geometry}

\usepackage{amscd}
\usepackage{amssymb}
\usepackage{amsmath}
\usepackage{amsfonts}
\usepackage{enumerate}
\usepackage[all,cmtip]{xy}
\usepackage{tikz-cd}
\usepackage{hyperref}
\usepackage{xcolor}
\usepackage{adjustbox}
\usepackage{graphicx}
\usepackage[T1]{fontenc}
\usepackage{lmodern}

\newtheorem{thm}{Theorem}[section]
\newtheorem{prop}[thm]{Proposition}
\newtheorem{lem}[thm]{Lemma}
\newtheorem{cor}[thm]{Corollary}

\theoremstyle{definition}

\newtheorem{exmp}[thm]{Example}
\newtheorem{que}[thm]{Question}
\newtheorem{rem}[thm]{Remark}

\newcommand{\Ss}{\mathbb{S}}
\newcommand{\Zz}{\mathbb{Z}}
\newcommand{\Zzp}{\mathbb{Z}_{(p)}}
\newcommand{\Ee}{\mathbb{E}}
\newcommand{\Ff}{\mathbb{F}_p}
\newcommand{\THH}{\mathrm{THH}}
\newcommand{\MU}{\mathrm{MU}}

\makeatletter
\@namedef{subjclassname@2020}{2020 Mathematics Subject Classification}
\makeatother

\subjclass[2020]{Primary 55Q51; Secondary 55P43}
\keywords{Chromatic height, chromatic redshift, topological Hochschild homology}

\numberwithin{equation}{section}

\title[Heights of $(T^n\otimes \Ff)^{hT^n}$, $(T^n\otimes \Zz)^{hT^n}$, and $(T^n\otimes j_p)^{hT^n}$]{Some chromatic height computations for homotopy fixed points of iterated topological Hochschild homology}

\author[J. Zha]{Jiaxi Zha}
\address{Department of Mathematics, Nankai University, No.94 Weijin Road, Tianjin 300071, P. R. China}
\email{1093913699@qq.com}

\begin{document}

\begin{abstract}
For an $\Ee_{\infty}$-ring spectrum $R$ of chromatic height $k$, the trace method and the chromatic redshift imply that $T^n$-fixed point of $n$-fold iterated topological Hochschild homology of $R$ has chromatic height at most $k+n$. In this paper, we compute these heights for $H\Ff$, $H\Zz$ and the $p$-complete connective image-of-$J$ spectrum $j_p$ for all primes and $n\geq1$, which are exactly $n-1$, $n$ and $n+1$, respectively. Our method is based on properties of Wilson spaces and relative topological Hochschild homology over $\MU$. 
\end{abstract}

\maketitle

\section{Introduction}
For a prime $p$, let $K(n)$ be the height $n$ Morava $K$-theory spectrum at $p$ with $K(-1)=\Ss$. In \cite{Hahn16}, Haha proved that, for an $\Ee_{\infty}$-ring spectrum $R$, if $K(n)\otimes R$ equivalent to zero then so is $K(n+1)\otimes R$; see also \cite[Theorem 9.4]{BSY22}. Consequently, we can define the chromatic height of $R$ by
	\[\mathrm{height}(R):=\sup\left\{n\geq-1\ |\ K(n)\otimes R\neq0 \right\}.\]
The chromatic redshift theorem then implies that the algebraic $K$-theory raises height by one. On the other hand, there is a trace map from algebraic $K$-theory to negative topological cyclic homology $K(R)\to \THH(R)^{hT}$ by \cite[Theorem 1.12]{BGT13}. Therefore, it is natural to ask whether the height of $K(R)$ can already be detected in $\THH(R)^{hT}$. For $p\geq5$, this question has a positive answer for the Adams summand and the connective complex $K$-theory by the computations of Ausoni-Rognes; see \cite{Aus10} and \cite{AR02}. More recently, Hahn-Wilson \cite{HW22} showed that $K(n+1)\otimes\THH(\mathrm{BP}\langle n\rangle)^{hT}$ is nonzero via relative topological Hochschild homology over $\MU$, which inspired our approach.

Iterating topological Hochschild homology gives a higher-dimensional version of this question. For an $\Ee_{\infty}$-ring spectrum $R$, let $T^n\otimes R$ denote the $n$-fold iterated topological Hochschild homology of $R$, 
\[\underbrace{\THH\circ\THH\circ\cdots\circ\THH}_{n}(R),\]
see Section \ref{sec2} for an explanation of this notation. Indeed, iterating the trace map, we obtain an $\Ee_{\infty}$-map $K^{(n)}(R)\longrightarrow(T^n\otimes R)^{hT^n},$ where $K^{(n)}$ denotes the $n$-fold iterated algebraic $K$-theory. Since $K(k)$ is an $\Ee_1$-ring spectrum for $k\geq-1$, the chromatic redshift theorem implies that $\mathrm{height}((T^n\otimes R)^{hT^n})\leq\mathrm{height}(R)+n$. Therefore, the main problem is to establish a lower bound for $\mathrm{height}((T^n\otimes R)^{hT^n})$. The main goal of this paper is to prove that, for every prime and $n\geq1$, the upper and lower bounds coincide for $R=H\Ff$, $H\Zz$ and $j_p$.

Our work is motivated by Veen's calculations in \cite{Vee18}. Veen proved that, for $p\geq5$ and $0\leq n\leq p-1$ or $p=3$ and $0\leq n\leq1$, the unit map
\[k(n)_*\longrightarrow k(n)_*(T^{n+1}\otimes H\Ff)^{hT^{n+1}}\]
sends $v_n$ to a nonzero element. Here $k(n)$ is the connective Morava $K$-theory spectrum. Recently, Fang \cite{Fang26} proved that, for $p\geq 5$ and $1\leq n\leq p-1$, the unit map
\[k(n)_*\longrightarrow k(n)_*(T^{n-1}\otimes j_p)^{hT^{n-1}}\]
also sends $v_n$ to a nonzero element. However, they did not prove that every positive power of $v_n$ has nonzero image under the corresponding unit map. Consequently, their computations do not suffice to determine the chromatic heights of these spectra in general. 

Veen's result \cite{Vee18}, which is used in  Fang's work \cite{Fang26}, relies on a detailed spectral sequence computation for $\pi_*(T^n\otimes H\Ff)$. Using Wilson spaces and relative topological Hochschild homology over $\MU$, we establish the required lower bounds without carrying out these computations. Our first main theorem is the following.

\begin{thm}\label{thm:main}
	For $n\geq 1$, the height of $(T^n\otimes H\Ff)^{hT^n}$ is $n-1$.
\end{thm}

We now describe the construction of the lower bound. First, we construct a map of spectra
\begin{equation}\label{eq:l/mu}
	f_n\colon(T^n\otimes H\Ff)^{hT^n}\to (T^{n-1}\otimes_{\MU}H\Zzp)^{hT^{n-1}}_p
\end{equation}
such that $f_n\circ\eta_1\simeq\eta_2$, where $\eta_1$ and $\eta_2$ are units of $(T^n\otimes H\Ff)^{hT^n}$ and $(T^{n-1}\otimes_{\MU}H\Zzp)^{hT^{n-1}}_p$, respectively. Next, we use the Wilson space in \cite{HRW25} to construct an $\Ee_{\infty}$-$\MU$-algebra $W(n-1)$ with a $T^{n-1}$-action and a $T^{n-1}$-equivariant map of $\Ee_{\infty}$-$\MU$-algebras
\begin{equation}\label{eq:wilson}
	g_n\colon T^{n-1}\otimes_{\MU}H\Zzp\to W(n-1).
\end{equation}
Finally, we prove that $K(n-1)\otimes W(n-1)^{hT^{n-1}}$ is nonzero, which implies that $K(n-1)\otimes(T^n\otimes H\Ff)^{hT^n}$ is also nonzero.

The construction of $W(n)$ directly gives the required lower bound for $H\mathbb Z$. Together with the upper bound above, this yields our second main theorem.

\begin{thm}\label{thm:main2}
	For $n\geq1$, the height of $(T^n\otimes H\Zz)^{hT^n}$ is $n$.
\end{thm}

Let $j_p$ be $p$-complete connective the image-of-$J$ spectrum. Using the same construction, with $\Zzp$ replaced by $\ell$, we can also obtain the following reslut.

\begin{thm}\label{thm:main3}
	For $n\geq1$, the height of $(T^n\otimes j_p)^{hT^n}$ is $n+1$. 
\end{thm}

Finally, as recalled above, Hahn-Wilson \cite{HW22} proved that $K(n+1)\otimes\THH(\mathrm{BP}\langle n\rangle)^{hT}$ is nonzero. The $\Ee_3$-ring structure which they constructed in \cite{HW22} allows $\THH$ to be iterated up to three times. This leads to the following question. Our method does not directly apply to this question, since our proof relies on \cite[Proposition 3.1.3]{HRW25}, which requires $\Ee_{\infty}$-ring structures.

\begin{que}
	Are $K(n+2)\otimes(\int_{T^2}\mathrm{BP}\langle n\rangle)^{hT^2}$ and $K(n+3)\otimes(\int_{T^3}\mathrm{BP}\langle n\rangle)^{hT^3}$ nonzero? Here $\int_{-}-$ is the factorization homology; see \cite[Appendix A]{HW22}.
\end{que}

\subsection*{Organization} In section \ref{sec2}, we recall the definition of iterated topological Hochschild homology and construct the map \eqref{eq:l/mu}. Then we construct the $\Ee_{\infty}$-ring spectrum $W(n)$ and the $T^n$-equivariant map of $\Ee_{\infty}$-ring spectra \eqref{eq:wilson} in Section \ref{sec3}. Then, in Section \ref{sec4}, we prove that $W(n)^{hT^n}$ is not $K(n)$-acyclic for $n\geq1$. Finally, we prove three main theorems in Section \ref{sec5}.
\subsection*{Conventions}
We work with $\infty$-categories throughout, following \cite{HTT,HA}. Therefore, all categories are understood to be $\infty$-categories unless otherwise specified. Moreover,
\begin{enumerate}
\item We write $\mathcal{S}$ for the category of spaces, and $\mathrm{Sp}$ for the category of spectra. The $\Ee_k$-ring spectra are called $\Ee_k$-rings. The maps between $\Ee_k$-ring are called $\Ee_k$-maps.
\item For an ordinal commutative ring $R$, the Eilenberg-MacLane spectrum $HR$ is also denoted by $R$. Moreover, $\MU$ denotes the $p$-localized complex cobordism spectrum.
\item For a spectrum $X$, we denote $X_p$ by its $p$-completion.
\item We write $\Lambda_R(x)$ and $\Gamma_R(x)$ for the exterior and divided power algebras over a commutative ring $R$, respectively. When $R=\Ff$, we simply write them as $\Lambda(x)$ and $\Gamma(x)$. 
\item All filtrations are assumed to be multiplicative.
\end{enumerate}

\section{The map $(T^{n+1}\otimes\Ff)^{hT^{n+1}}\to (T^n\otimes_{\MU}\Zzp)^{hT^n}_p$}\label{sec2}
In this section we construct the map \eqref{eq:l/mu}. First, we recall their definitions. By \cite[Corollary 4.4.4.9]{HTT} and \cite[Corollary 3.2.3.5]{HA}, there is a functor
\[\otimes_{\mathcal{C}}\colon \mathcal{S}\times\mathrm{CAlg}(\mathcal{C})\longrightarrow \mathrm{CAlg}(\mathcal{C})\]
for a presentable symmetric monoidal categoty $\mathcal{C}$. This functor maps $(X,R)$ to $\mathrm{colim}R_X$, where $R_X$ is the constant functor $X\to \mathrm{CAlg}(\mathcal{C})$ with the value $R$. When $\mathcal{C}$ is $\mathrm{Sp}$ or $\mathrm{Mod}_{\MU}$, we write $\otimes$ for $\otimes_{\mathrm{Sp}}$ and $\otimes_{\MU}$ for $\otimes_{\mathrm{Mod}_{\MU}}$, respectively. For a fixed $\Ee_{\infty}$-ring $R$ in $\mathcal{C}$, $-\otimes_{\mathcal{C}}R\colon\mathcal{S}\to\mathrm{CAlg}(\mathcal{C})$ preserves colimits by \cite[Example 2.5]{HY17}. Therefore, we deduce that $T\otimes\Ff\simeq\THH(R)$ and
\[T^n\otimes R\simeq (T\times T^{n-1})\otimes R\simeq T\otimes(T^{n-1}\otimes R)\simeq\THH(T^{n-1}\otimes R).\]
The second equivalence follows from the fact that $T\times T^{n-1}$ is the colimit of the constant functor $T\to \mathcal{S}$ with the value $T^{n-1}$. In other words, $T^n\otimes R$ is the $n$-fold iterated topological Hochschild homology of $R$. Moreover, we can consider $T^n$ as a functor $BT^n\to\mathcal{S}$. Compositing with $-\otimes R$, we obtain a $T^n$-action on $T^n\otimes R$.

Since the forgetful functor $\mathrm{Mod}_{\MU}\to\mathrm{Sp}$ is lax symmetric monoidal, it induces a functor $\mathrm{CAlg}(\mathrm{Mod}_{\MU})\to\mathrm{CAlg}(\mathrm{Sp})$. Hence for $X\in\mathcal{S}$ and $R\in\mathrm{CAlg}(\mathrm{Mod}_{\MU})$, we have a canonical comparison map
\begin{equation}\label{eq:forget}
	X\otimes R\rightarrow X\otimes_{\MU}R
\end{equation} 
of $\Ee_{\infty}$-rings, where we consider $R$ on the left as an $\Ee_{\infty}$-ring via the forgetful functor. 

Now, we construct the map \eqref{eq:l/mu}. For this, we need the following lemma.

\begin{lem}\label{lem:eta}
	Let $M$ be an $\MU$-module with trivial $T$-action. Then there is a map of $\MU$-modules $\epsilon\colon M^{tT}\to M$ such that the composition
	\[M\xrightarrow{c} M^{hT}\to M^{tT}\xrightarrow{\epsilon} M\]
	is the identity, where the first map is induced by the collapsing map $BT\to *$ and the second map is the canonical map from the fixed point to the Tate construction.
\end{lem}
\begin{proof}
	By \cite[Construction 2.4.4]{Rak26}, there is a cofiber sequence
	\[\Sigma M_{hT}\xrightarrow{\mathrm{Nm_M}} M^{hT}\to M^{tT}\]
	in the category of $\MU$-modules. Let $e_M\colon M^{hT}\to M$ be the map induced by the inclusion of the base point. Since $e_M\circ c$ is equivalent to the identity, it suffices to show that $e_M\circ\mathrm{Nm}_M$ is equivalent to zero in the category of $\MU$-modules. For this, we claim that there is a commutative diagram
	\[\xymatrix@C=6em@R=3em{
	M\otimes_{\MU}\Sigma\MU_{hT} \ar[r]^{M\otimes_{\MU}\mathrm{Nm}_{\MU}} \ar[d]_{\simeq} & M\otimes_{\MU} \MU^{hT} \ar[d] \\
	\Sigma M_{hT} \ar[r]_{\mathrm{Nm}_{M}} & M^{hT}.}\]
	Indeed, \cite[Construction 2.4.4]{Rak26} gives a commutative diagram 
	\[\xymatrix@C=6em@R=3em{
	M\otimes_{\MU}\Sigma\MU_{hT} \ar[r] \ar[d]_{\simeq} & M\otimes_{\MU} \MU \ar[d] \\
	\Sigma M_{hT} \ar[r] & M.}\]
	Now, the upper-right composite $M\otimes_{\MU}\Sigma\MU_{hT}\to M$ factors as 
	\[M\otimes_{\MU}\Sigma\MU_{hT}\to M\otimes_{\MU}\MU^{hT}\to M\otimes_{\MU}\MU\simeq M.\]
	Taking its adjunction implies the claim. Hence it suffices to show that $e_{\MU}\circ\mathrm{Nm}_{\MU}$ is equivalent to zero in the category of $\MU$-modules. This follows from the isomorphism
	\[[\Sigma\MU_{hT},\MU]_{\MU}\cong[\Sigma\MU\otimes\Sigma^{\infty}_+BT,\MU]_{\MU}\cong[\Sigma\Sigma^{\infty}_+BT,\MU]\cong0,\]
	since $\MU$ is complex orientable.
\end{proof}

\begin{prop}\label{prop:e0}
	For $n\geq0$, there exists an $\Ee_0$-map
	\[f_n\colon(T^{n+1}\otimes\Ff)^{hT^{n+1}}\to (T^n\otimes_{\MU}\Zzp)^{hT^n}_p.\]
\end{prop}
\begin{proof}
	Since $T\otimes\Ff$ is a cyclotomic spectrum, there is a canonical $T$-equivariant $\Ee_{\infty}$-map $\varphi\colon T\otimes\Ff\to (T\otimes\Ff)^{tC_p}$. By \cite[Corollary IV.4.13]{NS18}, there is a $T$-equivariant equivalence of $\Ee_{\infty}$-rings $(T\otimes\Ff)^{tC_p}\xrightarrow{\simeq}\Zz_p^{tC_p}.$ Compositing these two maps and applying $T^n\otimes-$, we obtain a $T^{n+1}$-equivariant $\Ee_{\infty}$-map
	\begin{equation}\label{eq:Zp}
		T^{n+1}\otimes\Ff\longrightarrow T^n\otimes\Zz_p^{tC_p}.
	\end{equation}
	Since the Tate construction is lax symmetric monoidal by \cite[Corollary I.3.9]{NS18}, there is a canonical $T^2$-equivariant $\Ee_{\infty}$-map $T\otimes R^{tC_p}\to (T\otimes R)^{tC_p}$ for $R\in\mathrm{Fun}(BT,\mathrm{CAlg}(\mathrm{Sp}))$, where the Tate construction on $T\otimes R$ is taken with respect to the $T$-action on $R$. Iterating this $n$ times to the right hand of the map \eqref{eq:Zp} and taking $T^{n+1}$-fixed point gives an $\Ee_{\infty}$-map
	\begin{equation*}
	\begin{split}
		(T^{n+1}\otimes\Ff)^{hT^{n+1}}\longrightarrow & ((T^n\otimes\Zz_p)^{tC_p})^{h\left(T^n\times(T/C_p)\right)} \\
		& \simeq \left((T^n\otimes\Zz_p)^{tT}_p\right)^{hT^n}\simeq \left(\left(T^n\otimes\Zzp\right)^{tT}\right)^{hT^n}_p,
	\end{split}
	\end{equation*}
	where the Tate construction is taken with respect to the $T$-action on $\Zz_p$ and $\Zzp$. The first equivalence follows from \cite[Lemma II.4.2]{NS18} and the second equivalence follows from the fact that $\Ss/p$ is finite. Now, the comparison map \eqref{eq:forget} and Lemma \ref{lem:eta} gives the desired map
	\begin{equation*}
		\begin{split}
			(T^{n+1}\otimes\Ff)^{hT^{n+1}}\longrightarrow\left(\left(T^n\otimes\Zzp\right)^{tT}\right)^{hT^n}_p & \longrightarrow\left(\left(T^n\otimes_{\MU}\Zzp\right)^{tT}\right)^{hT^n}_p \\
			& \xrightarrow{\ \ \epsilon\ \ }\left(T^n\otimes\Zzp\right)^{hT^n}_p.
		\end{split}
	\end{equation*}
	The map labelled by $\epsilon$ preserve the unit by Lemma \ref{lem:eta}, since it admits an $\Ee_{\infty}$-section. The other two maps also preserve the unit, since they are $\Ee_{\infty}$-maps. 
\end{proof}

Now, choose generators of $\MU_*$ such that $\MU_*\cong\Zzp[v_i,x_j\ |\ i\geq1,\ j\neq p^k-1]$ as in \cite[Lemma 2.0.5]{HW22}. Then the unit map $\eta\colon\MU\to\Zzp$ maps all generators to zero on homotopy groups. In particular, we can consider $v_i$ as an element of $\pi_{2p^i-1}(\Zzp/\MU)$ for $i\geq1$, where $\Zzp/\MU$ is the cofiber of $\eta$. Using the suspension map in \cite[Example A.2.4]{HW22}, we obtain the element $z_i=\sigma^2 v_i\in\pi_{2p^i}(T\otimes_{\MU}\Zzp)$ for $i\geq1$. 

We end this section by considering the relations among the images of $z_i$ in $\pi_*((T\otimes_{\MU}\Zzp)/p)$.

\begin{lem}\label{lem:modpv1}
	Let $\bar{z}_i$ be the image of $z_i$ in $\pi_*((T\otimes_{\MU}\Zzp)/p)$ for $i\geq1$. Then $\bar{z}_i\doteq\bar{z}_1^{p^{i-1}}$ for $i\geq1$. In particular, $\pi_*((T\otimes_{\MU}\Zzp)/p))$ is a free $\Ff[\bar{z}_1]$-module.
\end{lem}
\begin{proof}
	By base change along the canonical $\Ee_{\infty}$-map $\Zzp\to\Ff$, we obtain an $\Ee_{\infty}$-map $f\colon(T\otimes_{\MU}\Zzp)\otimes_{\Zzp}\Ff\to T\otimes_{\MU}\Ff$. The proof of Proposition 2.5.3 in \cite{HW22} implies that $f$ induces an injective map on homotopy groups. More precisely, it induces an injective map on the $E^{\infty}$-pages of Tor spectral sequences
	\[\Gamma(\bar{z}_i,\overline{\sigma^2x^{\prime}}_j\ |\ i\geq1,\ j\neq p^k-1)\ \longrightarrow \Gamma(\sigma^2v_0, z_i,\sigma^2x^{\prime}_j\ |\ i\geq1,\ j\neq p^k-1),\]
	which maps $\bar{z}_i$ to $z_i$ and $\overline{\sigma^2x^{\prime}}_j$ to $\sigma^2x^{\prime}_j$. Here $v_0=p\in\MU_*.$ Hence it suffices to show that $z_i\doteq z_2^{p^{i-2}}$ in $\pi_*(T\otimes_{\MU}\Ff)$ for $i\geq1$. By \cite[Proposition 2.4.2]{HW22}, we see that $\sigma^2v_i\doteq(\sigma^2v_0)^{p^i}$ in the $E^{\infty}$-page of the Tor spectral sequence for $\pi_*(T\otimes_{\MU}\Ff)$. Since these two elements have filtration degree one and $F_0\pi_{2p^i}(T\otimes_{\MU}\Ff)=\pi_{2p^i}(\Ff)=0$, we obtain that $\sigma^2v_i\doteq(\sigma^2v_0)^{p^i}$ in $\pi_*(T\otimes_{\MU}\Ff)$. Hence $z_i\doteq((\sigma^2v_0)^{p})^{p^{i-1}}\doteq z_1^{p^{i-1}}$. In particular, $\bar{z}_1^m$ is nonzero for all $m\geq0$, and the freeness result follows immediately.
\end{proof}

\section{The construction of $W(n)$ and the map \eqref{eq:wilson}}\label{sec3}
In this section, we use the Wilson space to construct the $\Ee_{\infty}$-ring $W(n)$ mentioned in the introduction and study some of its properties.

 Recall that, for $i\geq1$, the $2i$-th Wilson space is defined to be $W_{2i}:=\Omega^{\infty}\Sigma^{2i}\MU$; see \cite{Wil73}. By \cite[Theorem 3.1.2]{HRW25}, $\MU_*(W_{2i})$ is a polynomial $\MU_*$-algebra, with polynomial generators in even degrees. Moreover, since $\MU$ is a connective spectrum, the K\"unneth formula implies that $\MU_{2i}(W_{2i})\cong\Zzp$ and $\MU_k(W_{2i})=0$ for $k<2i$. For a connective even $\Ee_{\infty}$-$\MU$-algebra $R$ and $x\in\pi_{2i}(R)$, \cite[Proposition 3.1.3]{HRW25} and its proof imply that there is an $\Ee_{\infty}$-map
\begin{equation}\label{eq:w2itor}
	\MU\otimes\Sigma^{\infty}_+W_{2i}\longrightarrow R,
\end{equation}
which maps the generator in $\MU_{2i}(W_{2i})\cong\Zzp$ to $x\in\pi_{2i}(R)$. Using this property, we construct $W(n)$ by an induction on $n$. For this, we need the following lemma.

\begin{lem}\label{lem:ind}
	Let $R$ be a connective even $\Ee_{\infty}$-$\MU$-algebra equipped with a $T^n$-action and an $\Ee_{\infty}$-$\MU$-algebra map $\Zzp\to R$. Suppose that
	\begin{enumerate}
		\item[$(1)$] $\pi_*(R)$ is isomorphic to $\Zzp[x_{\alpha}]/(f_{\beta})$ as $\Zzp$-algebras with $|x_{\alpha}|,|f_{\beta}|>0$,
		\item[$(2)$] for any $k>0$, there are only finitely many $x_{\alpha}$ and $f_{\beta}$ of degree $2k$,
		\item[$(3)$] every finite subfamily of $\{f_{\beta}\}$ is a regular sequence.
	\end{enumerate}
	Then there is a $T^n$-equivariant $\Ee_{\infty}$-$\MU$-algebra map $\MU\otimes\Sigma^{\infty}_+W\to R$, where $W$ is a product of Wilson spaces $W_{2i}$ with finitely many factors for each $i\geq1$ and $\MU\otimes\Sigma^{\infty}_+W$ is equipped with the trivial $T^n$-action. The induced map on homotopy groups is surjective. Its kernel, denoted by $I$, is generated by elements of positive even degree, with only finitely many in each degree, such that every finite subfamily of generators is a regular sequence in $\MU_*(W)$. Moreover, $I/I^2$ is a free $\pi_*(R)$-module with a basis containing $v_i$ for all $i\geq1$.
\end{lem}
\begin{proof}
	Since $R$ is even, the homotopy fixed point spectral sequence for $\pi_*(R^{hT^n})$ collapses. Hence $R^{hT^n}$ is also even and the canonical map $\pi_*(R^{hT^n})\to \pi_*(R)$ is surjective. Consequently, considering $x_{\alpha}$ as an element of $\pi_*(R)$, we can lift it to $\pi_*(R^{hT^n})$. By the map \eqref{eq:w2itor}, we obtain an $\Ee_{\infty}$-map $\MU\otimes\Sigma^{\infty}_+W_{|x_{\alpha}|}\to R^{hT^n}$ for each element $x_{\alpha}\in\pi_*(R)$. Set $W=\left(\prod_{\alpha}W_{|x_{\alpha}|}\right)$. Since there are only finitely many $x_{\alpha}$ in each positive degree, there are only finitely many $W_{2i}$-factors for each $i\geq1$. Multiplying these maps together in the category of $\MU$-algebras and using the adjunction between the trivial action and the homotopy fixed point give the desired $T^n$-equivariant $\Ee_{\infty}$-map 
	\[g\colon\MU\otimes\Sigma^{\infty}_+W\longrightarrow R.\]
	By \cite[Theorem 3.1.2]{HRW25}, $\MU_*(W)\cong\Zzp[v_i,x^{\prime}_j,a_{\alpha},c_{\gamma}\ |\ i\geq1,\ j\neq2p^k-1]$ and $g_*(a_{\alpha})=x_{\alpha}$ by construction. In particular, $g_*$ is surjective. For each $\gamma$, choose a polynomial $F_{\gamma}(x_{\alpha})\in\Zzp[x_{\alpha}]$ representing $g_*(c_{\gamma})\in\pi_*(R)\cong\Zzp[x_{\alpha}]/(f_{\beta})$ and let $F_{\gamma}(a_{\alpha})$ be the polynomial obtained by replacing $x_{\alpha}$ by $a_{\alpha}$. Then we have that $I=\left(v_i,x^{\prime}_j,f_{\beta}(a_{\alpha}),c_{\gamma}-F_{\gamma}(a_{\alpha})\right)$. 
	
	Let $\{h_{\delta}\}=\{v_i,x^{\prime}_j,f_{\beta}(a_{\alpha}),c_{\gamma}-F_{\gamma}(a_{\alpha})\}$, and consider the $\MU_*(W)$-linear map
	\[\bigoplus_{\delta}\MU_*(W)\{h_{\delta}\}\longrightarrow I/I^2,\]
	which maps $ah_{\delta}$ to $\overline{ah}_{\delta}$ for each $\delta$ and $a\in \MU_*(W)$. If $a\in I$, then $ah_{\delta}\in I/I^2$. Hence it descends to a $\pi_*(R)$-linear map
	\begin{equation}\label{eq:I/I2}
	\bigoplus_{\delta}\pi_*(R)\{h_{\delta}\}\longrightarrow I/I^2.		
	\end{equation}
	Since every finite subfamily of $\{f_{\beta}\}$ is a regular sequence, so is that for $\{h_{\delta}\}$.

	This map is surjective, since every element in $I/I^2$ is a $\pi_*(R)$-linear combination of $\overline{h}_{\delta}$. Next, suppose that $\Sigma_i a_ih_i\in I^2$ with $a_i\in\MU_*(W)$ and $h_i\in\{h_{\delta}\}$. Since the sum is finite, we can choose $k$ large enough such that
	\[\sum_{1\leq i\leq k} a_ih_i=\sum_{1\leq j\leq m}u_jw_j\]
	with $u_j,w_j\in(h_1,\cdots,h_k)$ for all $1\leq j\leq m$, by setting $a_i=0$ for newly added generators. Hence we deduce that
	\[\sum_{1\leq i\leq k} a_ih_i=\sum_{1\leq j\leq m}u_jw_j=\sum_{1\leq j\leq m}\left(\sum_{1\leq i\leq k}c_{i,j}h_i\right)w_j=\sum_{1\leq i\leq k}\left(\sum_{1\leq j \leq m}c_{i,j}w_j\right)h_i,\]
	for some $c_{i,j}\in\MU_*(W)$. That is, $\sum_{1\leq i\leq k} (a_i-b_i)h_i=0$ with $b_i\in(h_1,\cdots,h_k)$. Since $\{h_1,\cdots,h_k\}$ is a regular sequence, we obtain that $a_i-b_i\in(h_1,\cdots,h_k)$ for $1\leq i\leq k$. In particular, $a_i\in I$, which implies that the map \eqref{eq:I/I2} is an isomorphism. Since $v_i\in \{h_{\delta}\}$ for all $i\geq1$, the result follows.
\end{proof}

\begin{prop}\label{prop:ind}
	Let $R$ be an $\Ee_{\infty}$-$\MU$-algebra as in Lemma \ref{lem:ind}. Then $T\otimes_{\MU\otimes\Sigma^{\infty}_+W}R$ satisfies the hypotheses of Lemma \ref{lem:ind}, with $n$ replaced by $n+1$. Moreover, $\pi_*(T\otimes_{\MU\otimes\Sigma^{\infty}_+W}R)$ admits a filtration, whose associated graded is
	\[\mathrm{gr}(\pi_*(T\otimes_{\MU\otimes\Sigma^{\infty}_+W}R))=\Gamma_{\pi_*(R)}(\sigma^2h_{\delta}\ |\ h_{\delta}\in\{h_{\delta}\}),\]
	where $\{h_{\delta}\}$ is the set given in the proof of Lemma \ref{lem:ind}. 
\end{prop}
\begin{proof}
	Since any finite subfamily of $\{h_{\delta}\}$ is a regular sequence, \cite[Proposition 3.6]{Ang08} implies that there is an isomorphism of $\pi_*(R)$-algebras 
	\[\pi_*(R\otimes_{\MU\otimes\Sigma^{\infty}_+W}R)\cong\Lambda_{\pi_*(R)}(\sigma h_{\delta}).\]
	Therefore, the Tor spectral spectral for $\pi_*(T\otimes_{\MU\otimes\Sigma^{\infty}_+W}R)$ has the $E^2$-page
	\[\mathrm{Tor}^{\Lambda_{\pi_*(R)}(\sigma h_{\delta})}_{*,*}(\pi_*(R),\pi_*(R))\cong \Gamma_{\pi_*(R)}(\sigma^2h_{\delta}).\]
	Since the total degree of $\sigma^2h_{\delta}$ is even, this spectral sequence also collapses at the $E^2$-page, which gives the desired associated graded. In particular, $T\otimes_{\MU\otimes\Sigma^{\infty}_+W}R$ is connective even with homotopy groups free over $\pi_*(R)$. 
	
	For $a\geq0$, let $x_{\delta,a}=\gamma_{p^a}(\sigma^2h_{\delta})$. Since $p$ in nonzero in $\pi_*(R)$, we have that $x_{\delta,a}^p=\frac{(p^{a+1})!}{((p^a)!)^p}x_{\delta,a+1}=pu_ax_{\delta,a+1}$ with $u_a\in\Zzp^{\times}.$ Hence there is an isomorphism
	\[\Gamma_{\pi_*(R)}(\sigma^2h_{\delta})\cong\pi_*(R)[x_{\delta,a}]/(x_{\delta,a}^p-pu_ax_{\delta,a+1}),\]
	and we can choose lifts $\widetilde{x}_{\delta,a}$ of $x_{\delta,a}$ in $\pi_*(T\otimes_{\MU\otimes\Sigma^{\infty}_+W}R)$. Consider the map
	\[F\colon\pi_*(R)[X_{\delta,a}]\longrightarrow\pi_*(T\otimes_{\MU\otimes\Sigma^{\infty}_+W}R),\qquad X_{\delta,a}\mapsto \widetilde{x}_{\delta,a}\]
	and define a multiplicative filtration on $\pi_*(R)[X_{\delta,a}]$ by setting $\mathrm{fil}(\pi_*(R))=0$ and $\mathrm{fil}(X_{\delta,a})=p^a$ for $a\geq0$. Then $F$ preserves the filtration and $\mathrm{gr}F$ is surjective. Since the filtrations are finite in a fixed degree on both sides and $F_0(F)$ is an isomorphism, $F$ is surjective. Moreover, since $F(X_{\delta,a}^p-pu_aX_{\delta,a+1})=\widetilde{x}_{\delta,a}^p-pu_a\widetilde{x}_{\delta,a+1}\in F_{p^{a+1}-1}(\pi_*(T\otimes_{\MU\otimes\Sigma^{\infty}_+W}R))$, we can choose  $G\in F_{p^{a+1}-1}(\pi_*(R)[X_{\delta,a}])$ such that $F(G)=F(X_{\delta,a}^p-pu_aX_{\delta,a+1}),$ which implies that
	\[F(X_{\delta,a}^p-pu_aX_{\delta,a+1}-G)\in F_{p^{a+1}-2}(\pi_*(T\otimes_{\MU\otimes\Sigma^{\infty}_+W}R)).\]
	Iterating this $p^{a+1}$ times, we obtain a polynomial $G_{\delta,a}$ suhch that $F(X_{\delta,a}^p-pu_aX_{\delta,a+1}-G_{\delta,a})=0$. Hence $F$ descends to a map
	\[\widetilde{F}\colon\pi_*(R)[X_{\delta,a}]/(X_{\delta,a}^p-pu_aX_{\delta,a+1}-G_{\delta,a})\longrightarrow\pi_*(T\otimes_{\MU\otimes\Sigma^{\infty}_+W}R)\]
	with $\mathrm{fil}(G_{\delta,a})<p^{a+1}$. In particular, $\widetilde{F}$ induces an isomorphism on the associated graded, which implies that $\widetilde{F}$ is injective. Consequently, we obtain that
	\[\pi_*(T\otimes_{\MU\otimes\Sigma^{\infty}_+W}R)\cong\Zzp[x_{\alpha},X_{\delta,a}]/(f_{\beta},X_{\delta,a}^p-pu_aX_{\delta,a+1}-G_{\delta,a}).\]
	It remains to verify $(3)$ in Lemma \ref{lem:ind}, which can be checked in the associated graded by the following lemma.
\end{proof}

\begin{lem}\label{lem:gr and reg}
	Let $A$ be a commutative ring which admits an exhaustive filtration and $x_1,\cdots,x_k$ be a sequence of elements of $A$ with $x_i\in F_{d_i}A\setminus F_{d_i-1}A$. Let $\bar{x}_i\in \mathrm{gr}_{d_i}(A)$ be the image of $x_i$ in the associated graded of $A$. If $\bar{x}_1,\cdots,\bar{x}_k$ is a regular sequence in $\mathrm{gr}(A)$, then $x_1,\cdots,x_k$ is a regular sequence in $A$ and 
	\[\mathrm{gr}(A/(x_1,\cdots,x_k))\cong\mathrm{gr}(A)/(\bar{x}_1,\cdots,\bar{x}_k).\]
\end{lem}
\begin{proof}
	We first consider that case $k=1$. Let $y\in F_sA\setminus F_{s-1}A$ be a nonzero element. Then $\bar{x}_1\bar{y}$ is also nonzero is $\mathrm{gr}(A)$, since $\bar{x}_1$ is a nonzerodivisor. In particular, $x_1$ is a nonzerodivisor and $x_1y$ lies in $F_tA$ for some $t\geq d_1$ if and only if $s\leq t-d_1$. That is, $(x_1)\cap F_tA=x_1F_{t-d_1}A$. Let $p:A\to A/x_1$ be the canonical projection. Then we have that
	\[\mathrm{Ker}(\mathrm{gr}_s(p))\cong\frac{(x_1)\cap F_sA+F_{s-1}A}{F_{s-1}A} \cong\frac{x_1F_{s-d_1}A+F_{s-1}A}{F_{s-1}A}\cong \bar{x}_1\mathrm{gr}_{s-d_1}(A),\]
	which implies that $\mathrm{gr}(A/x_1)\cong\mathrm{gr}(A)/\bar{x}_1$. Consequently, repeating this argument and inducting on $k$ yield the desired result.
\end{proof}

For $1\leq i\leq n$, the inclusion of $i$-th $T$ into $T^n$ gives a map $\iota_i\colon T\otimes_{MU}\Zzp\to T^n\otimes_{\MU}\Zzp$. We denote the image of $z_j\in\pi_*(T\otimes_{\MU}\Zzp)$ under $(\iota_i)_*$ by $\zeta_{i,j}$ for $j\geq1$. For convenience, we also write $\zeta_i$ for $\zeta_{i,1}$ . Now, we can construct the $\Ee_{\infty}$-$\MU$-algebra $W(n)$ by induction on $n$. 

\begin{prop}\label{prop:wn}
	For $n\geq1$, there is an even $\Ee_{\infty}$-$\MU$-algebra $W(n)$ equipped with a $T^n$-action and a $T^n$-equivariant $\Ee_{\infty}$-$\MU$-algebra map $g_n\colon T^n\otimes_{\MU}\Zzp\to W(n)$. Moreover, $\pi_*(W(n))$ is a free $\Zzp$-module and $\pi_*(W(n)/p)$ is a free $\Ff[\zeta_1,\cdots,\zeta_n]$-module. Here, for $1\leq i\leq n$, we consider $\zeta_i$ as an element of $\pi_*(W(n)/p)$ via $(g_n)_*$.
\end{prop}
\begin{proof}
	For $n=1$, we take $W(1)=T\otimes_{\MU}\Zzp$ and $g_1=id$. The desired properties follows from \cite[Theorem 2.5.3]{HW22} and Lemma \ref{lem:modpv1}. Moreover, \cite[Theorem 2.5.3]{HW22} also implies that $W(1)$ satisfies the hypotheses of Lemma \ref{lem:ind}.
	
	Next, suppose that we have constructed $W(n)$ with the desired properties and satisfying the hypotheses of Lemma \ref{lem:ind}. Then we can choose a $T^n$-equivariant $\Ee_{\infty}$-$\MU$-algebra map $\MU\otimes\Sigma^{\infty}_+W\to W(n)$ as in Lemma \ref{lem:ind}. Using the inductive hypotheses, we take 
	\[W(n+1):=T\otimes_{\MU\otimes\Sigma^{\infty}_+W}W(n)\]
	and 
	\[g_{n+1}\colon T^{n+1}\otimes_{\MU}\Zzp\simeq T\otimes_{\MU}(T^n\otimes_{\MU}\Zzp)\xrightarrow{T\otimes_{\MU} g_n} T\otimes_{\MU}W(n)\xrightarrow{T\otimes_{\eta}W(n)}W(n+1),\]
	where $\eta$ is the unit map $\MU\to\MU\otimes\Sigma^{\infty}_+W$. By Proposition \ref{prop:ind}, $W(n+1)$ satisfies the hypotheses of Lemma \ref{lem:ind}, and there is a filtration on $\pi_*(W(n+1))$ such that
	\[\mathrm{gr}(\pi_*(W(n+1)))=\Gamma_{\pi_*(W(n))}(\sigma^2h_{\delta}),\]
	where $\{h_{\delta}\}$ is set of the generators of $\mathrm{Ker}(\MU_*(W)\to\pi_*(W(n))$ containing $v_i$ for $i\geq1$. Therefore, $\pi_*(W(n+1))$ is free over $\pi_*(W(n))$, hence over $\Zzp$. In particular, $p$ is a nonzerodivisor. Now, Lemma \ref{lem:gr and reg} implies that
	\[\mathrm{gr}(\pi_*(W(n+1)/p))\cong\Gamma_{\pi_*(W(n)/p)}(\sigma^2h_{\delta})\cong \pi_*(W(n)/p)[x_{\delta,k}\ |\ k\geq0]/(x_{\delta,k}^p),\]
	since $p=0$ in $\pi_*(W(n)/p)$. Here $x_{\delta,k}$ represents $\gamma_{p^k}(\sigma^2h_{\delta})$ for each $k\geq0$. Since $\{h_{\delta}\}$ contains $v_i$ for $i\geq1$, $\{x_{\delta,0}\}$ contains $\zeta_{n+1,i}$ for $i\geq1$. By Lemma \ref{lem:modpv1}, we have that $\zeta_{n+1,i}\doteq\zeta_{n+1}^{p^{i-1}}$. Let $B_{n+1}$ be the set of monomials in $\pi_*(W(n)/p)[x_{\delta,k}\ |\ k\geq0]/(x_{\delta,k}^p)$ not containing $\zeta_{n+1,i}$ for $i\geq1$. Choose lifts $\widetilde{b}\in\pi_*(W(n+1)/p)$ of $b\in B_n$ with $\widetilde{1}=1$ and define
	\[\Phi\colon\bigoplus_{b\in B_{n+1}}\pi_{*-|b|}(W(n)/p)[\zeta]\to \pi_*(W(n+1)/p)\]
	which maps $x\zeta^m$ to $x\zeta_{n+1}^m\widetilde{b}$ in the summand labelled by $b$ for $m\geq0$, where $x\in\pi_{*-|b|}(W(n)/p)$ and $|\zeta|=2p$. For $m=d_0+d_1p+\cdots d_lp^l$ with $0\leq d_0\cdots d_l<p$, let $\zeta^m$ have filtration $d_0+\cdots+d_l+\mathrm{fil}(b)$ and elements in $\pi_{*-|b|}(W(n)/p)$ have filtration degree zero  in a $b$-summand. Then $\Phi$ is a filtered map, and $\mathrm{gr}(\Phi)$ is an isomorphism by Lemma \ref{lem:modpv1}. Since both filtrations are finite in a fixed degree and $F_0\Phi$ is an isomorphism, $\Phi$ is also an isomorphism. In particular, we obtain an injective $\pi_*(W(n)/p)$-algebra map $\pi_*(W(n)/p)[\zeta]\to\pi_*(W(n+1)/p)$ which maps $\zeta^m$ to $\zeta_{n+1}^m$ for $m\geq0$. Hence $\pi_*(W(n+1)/p)$ is a free $\pi_*(W(n)/p)[\zeta_{n+1}]$-module. By inductive hypotheses, $\pi_*(W(n)/p)$ is a free $\Ff[\zeta_1,\cdots,\zeta_n]$-module. Consequently, $\pi_*(W(n+1)/p)$ is a free $\Ff[\zeta_1,\cdots,\zeta_{n+1}]$-module.
\end{proof}

\section{The chromatic height of $W(n)$}\label{sec4}
In this section we show that $K(n)\otimes W(n)^{hT^n}$ is nonzero for $n\geq0$ and finish the proof of Theorem \ref{thm:main}. We begin with two algebraic lemmas, which deals with regular sequences.

\begin{lem}\label{lem:reg in poly}
	Let $A$ be an $\Ff[a_1,\cdots,a_n]$-algebra whose underlying module is free. Then the sequence
	\[\left\{c_j=\sum_{i=1}^nb_ia_i^{p^{j-1}}\ |\ 1\leq j\leq n\right\}\]
	is a regular sequence in $A[b_1,\cdots,b_n]$.
\end{lem}
\begin{proof}
	Since $A$ is a free $\Ff[a_1,\cdots,a_n]$-module,\ it suffices to show that $\{c_j\}$ is regular sequence in $\Ff[a_1,\cdots,a_n,b_1\cdots,b_n]$. Set the weight of $a_i$ to be $2p^{i-1}$ and the weight of $b_i$ to be $p^{2n}-p^{2i-2}$ for $1\leq i\leq n$, and extend the weight additively to monomials. Moreover, choose a monomial order such that $w(M)<w(N)$ implies $M<N$. Then, for $1\leq i,j\leq n$, we have that
	\[w(b_ia_i^{p^{j-1}})=p^{2n}+p^{2j-2}-(p^{i-1}-p^{j-1})^2\]
	and the monomial in $c_j$ of maximal order is $b_ja_j^{p^{j-1}}$. Since $\{b_1a_1,\cdots,b_na_n^{p^{n-1}}\}$ is a regular sequence, so is $\{c_j\}$ by \cite[Proposition 15.15]{Eis95}
\end{proof}

The next lemma is analogous to Lemma \ref{lem:gr and reg}, with the increasing filtration replaced by a decreasing filtration. It is used to handle the homotopy fixed point spectral sequence.

\begin{lem}\label{lem:gr and reg2}
	Let $B$ be a commutative ring which admits a complete separated filtration and $y_1,\cdots,y_k$ be a sequence of elements of $B$ with $y_i\in F^{d_i}\setminus F^{d_i+1}$. Let $\bar{y}_i\in\mathrm{gr}_{d_i}(B)$ be the image of $y_i$ in the associated graded of $B$. If $\bar{y}_1,\cdots,\bar{y}_k$ is a regular sequence in $\mathrm{gr}(B)$, then $y_1,\cdots,y_k$ is a regular sequence in $B$ and
	\[\mathrm{gr}(B/(y_1,\cdots,y_k))\cong\mathrm{gr}(B)/(\bar{y}_1,\cdots,\bar{y}_k).\]
\end{lem}
\begin{proof}
	By separation, any nonzero element lies in $F^sB\setminus F^{s+1}B$ for some $s$. Hence the same argument as the proof of Lemma \ref{lem:gr and reg} implies that $y_1$ is a nonzerodivisor, $(y_1)\cap F^sB=y_1F^{s-d_1}B$ for $s\geq d_1$ and $\mathrm{gr}(B/(y_1))\cong\mathrm{gr}(B)/(\bar{y}_1)$. Let $z\in\bigcap_s((y_1)+F^sB)$, and choose $w_s\in B$ such that $z-y_1w_s\in F^sB$. Then $y_1(w_s-w_{s+1})\in F^sB$, which implies that $w_s-w_{s+1}\in F^{s-d_1}B$ for $s\geq d_1$. Hence $w_s$ converges to $w\in B$ by completeness, and $z=y_1w$ by separation, which implies that the induced filtration on $B/(y_1)$ is separated. Moreover, we have the following exact sequence
	\[0\longrightarrow B/F^{s-d_1}B\xrightarrow{\ \bar{y}_1\ }B/F^sB\longrightarrow B/((y_1)+F^sB)\longrightarrow0\]
	for $s\geq d_1$. Since $B/F^{s+1-d_1}\to B/F^{s-d_1}$ is surjective for all $s\geq d_1$, $\mathrm{lim}^1B/F^{s-d_1}B$ is zero. Hence the completeness of $B$ implies that $B/(y_1)\cong\mathrm{lim}B/((y_1)+F^sB)$. That is, $B/(y_1)$ is complete. Repeating this argument and inducting on $k$ yield the desired result. 
\end{proof}

\begin{prop}\label{prop:k(n)wn}
	For $n\geq1$, $K(n)\otimes W(n)^{hT^n}\not\simeq0$.
\end{prop}
\begin{proof}
	Since $W(n)$ is even by Proposition \ref{prop:wn}, the cellular structure on $BT^n$ induces a complete separated filtration on $\pi_*(W(n)^{hT^n})$ whose associated graded is 
	\[\mathrm{gr}(\pi_*(W(n)^{hT^n})\cong\pi_*(W(n))[t_1,\cdots,t_n],\]
	where $t_i$ has filtration degree on and total degree $-2$ for $1\leq i\leq n$. By Proposition \ref{prop:wn}, $\pi_*(W(n))$ is a free $\Zz_{(p)}$-module. Hence $p$ is a nonzerodivisor in the associated graded. Lemma \ref{lem:gr and reg2} implies that $p$ is a nonzerodivisor in $\pi_*(W(n)^{hT^n})$ and the induced filtration on $\pi_*(W(n)^{hT^n}/p)$ has the associated graded
	\[\mathrm{gr}(\pi_*(W(n)^{hT^n}/p)\cong\pi_*(W(n)/p)[t_1,\cdots,t_n].\]
	Now, we have the following commutative diagram
	\[\xymatrix{
	\MU \ar[r] \ar[d] & W(n)^{hT^n} \ar[d] \\
	\Zz_{(p)} \ar[r] & W(n)}\]
	by Proposition \ref{prop:wn}. Hence the image of $v_i\in\MU_*$ lies in $F^1\pi_*(W(n)^{hT^n})=\mathrm{Ker}(\pi_*(W(n)^{hT^n})\to\pi_*(W(n)))$ for $i\geq1$. By \cite[Lemma A.4.1]{HW22}, there is a commutative diagram for each $1\leq i\leq n$
	\[\xymatrix{
	\Sigma^{-1}(\Zz_{(p)}/\MU) \ar[d]_{\sigma^2} \ar[r] & \MU \ar[r] \ar[d] & W(n) \ar[d] \\
	\Sigma^{-2}T\otimes_{\MU}\Zz_{(p)} \ar[r]^{t} & \mathrm{lim}_{\mathbb{C}P^1}T\otimes_{\MU}\Zz_{(p)} \ar[r] & W(n)^{hT_i},}\]
	where we write $T_i$ for $i$th factor in $T^n$. Therefore, using Lemma \ref{lem:ind} and Proposition \ref{prop:ind}, we see that $v_j\in \pi_{2p^j-2}(\MU)$ is sent to $t_i\zeta_{i,j}\in\mathrm{gr}(\pi_*(W(n)^{hT_i}))$ for $1\leq i,j\leq n$, which implies that
	\[v_j=\sum_{i=1}^nt_i\zeta_{i,j}\doteq\sum_{i=1}^nt_i\zeta_i^{p^{j-1}}\]
	in $\mathrm{gr}(\pi_*(W(n)^{hT^n}/p))$ by Lemma \ref{lem:modpv1}. Here we also wirte $v_j$ for its image in $\pi_*(W(n)^{hT^n})$. Since $\pi_*(W(n)/p)$ is a free $\Ff[\zeta_1,\cdots,\zeta_n]$-module by Proposition \ref{prop:wn}, Lemma \ref{lem:reg in poly} and Lemma \ref{lem:gr and reg2} implies that $v_1,\cdots,v_{n-1}$ is a regular sequence in $\pi_*(W(n)^{hT^n}/p)$. Since $W(n)^{hT^n}/p$ is even and the degrees of $v_1,\cdots,v_{n-1}$ are even, we obtain that
	\[\pi_*(W(n)^{hT^n}/(p,v_1,\cdots,v_{n-1}))\cong\pi_*(W(n)^{hT^n})/(p,v_1,\cdots,v_{n-1}).\]
	Lemma \ref{lem:reg in poly} and Lemma \ref{lem:gr and reg2} also implies that $v_n$ is a nonzerodivisor in 
	\[\pi_*(W(n)^{hT^n}/(p,v_1,\cdots,v_{n-1})).\]
	Hence $W(n)^{hT^n}/(p,v_1,\cdots,v_{n-1})[v_n^{-1}]$ is nonzero. Now the result follows from the next lemma.
\end{proof}

\begin{lem}
	Let $M$ be a $\MU$-module with $M/(p,v_1,\cdots,v_{n-1})[v_n^{-1}]$ nozero. Then $K(n)\otimes M\not\simeq0.$
\end{lem}
\begin{proof}
	For $1\leq i\leq n-1$, we have the following commutaive diagram in the category of $\MU$-modules
	\[\xymatrix{
	\Sigma^{4p^i-4}\MU \ar[d]^{v_i} \ar[r]^{v_i} & \Sigma^{2p^i-2}\MU \ar[r]^{p} \ar[d]^{v_i} & \Sigma^{2p^i-2}\MU/v_i \ar[d]^{v_i} \\
	\Sigma^{2p^2-2}\MU \ar[r]^{v_i} & \MU \ar[r]^{p} & \MU/v_i.}\]
	Since $v_i\circ p\simeq p\circ v_i\simeq0$, the map $v_i\colon\Sigma^{2p^i-2}\MU/v_i\to\MU/v_i$ factors through $\Sigma^{4p^i-3}\MU\to\MU/v_i$. Hence $v_i$ acts trivially on $\MU/v_i$, as $\pi_{4p^i-3}(\MU/v_i)=0$. Using the equivalence 
	\[M/(p,v_1,\cdots,v_{n-1})\simeq M\otimes_{\MU}(\MU/p)\otimes_{\MU}(\MU/v_1)\otimes_{\MU}\cdots\otimes_{\MU}(\MU/v_{n-1}),\]
	we see that $v_i$ acts trivially on $M/(p,v_1,\cdots,v_{n-1})$ for all $1\leq i\leq n-1$.
	
	Next, since $\MU/(p,v_1,\cdots,v_{n-1})$ is dualizable in the category of $\MU$-modules, we obtain that
	\begin{equation*}
		\begin{split}
			\MU/(p,v_1,\cdots,v_{n-1})\otimes_{\MU} L_{K(n)}M & \simeq L_{K(n)}(M/(p,v_1,\cdots,v_{n-1}))\\
			& \simeq (M/(p,v_1,\cdots,v_{n-1})[v_n^{-1}])^{\wedge}_{I_n}
		\end{split}
	\end{equation*}
	by \cite[Proposition 3.7]{BS16}, where $(-)^{\wedge}_{I_n}$ is the completion with respect to $(p,v_1,\cdots,v_{n-1})$. Since $v_i$ acts trivially on $M/(p,v_1,\cdots,v_{n-1})$ for all $1\leq i\leq n-1$, we obtain that
	\[(M/(p,v_1,\cdots,v_{n-1})[v_n^{-1}])^{\wedge}_{I_n}\simeq M/(p,v_1,\cdots,v_{n-1})[v_n^{-1}]\not\simeq0.\]
	Hence $L_{K(n)}M\not\simeq0$, which implies that $K(n)\otimes M\not\simeq0$.
\end{proof}

Since multiplying by $p$ is an equivalence in the fiber of $X\to X_p$, $K(n)\otimes X_p\simeq K(n)\otimes X$. Hence we obtain the following corollary.

\begin{cor}\label{cor:wnp}
	For $n\geq1$, $K(n)\otimes W(n)^{hT^n}_p\not\simeq0.$
\end{cor}

\section{Proofs of the main theorems}\label{sec5}
We begin with the proofs of Theorem \ref{thm:main} and Theorem \ref{thm:main2}.
\begin{proof}[Proof of Theorem \ref{thm:main}]
	By \cite[Theorem 3.3]{Fang26}, the height of $(T^n\otimes\Ff)^{hT^n}$ is at most $n-1$. It remains to show that $K(n-1)\otimes (T^n\otimes\Ff)^{hT^n}\not\simeq0$. For the case $n=1$, there is an isomorphism $\pi_0((T\otimes\Ff)^{hT})\cong\Zz_p$ by \cite[Proposition IV.4.6]{NS18}. In particular $\mathbb{Q}\otimes(T\otimes\Ff)^{hT}$ is nonzero. 
	
	Next, assume that $n\geq2$. By Proposition \ref{prop:e0}, there is an $\Ee_0$-map $(T^n\otimes\Ff)^{hT^n}\to (T^{n-1}\otimes_{\MU}\Zzp)^{hT^{n-1}}_p$. Compositing with the $\Ee_{\infty}$-map $(T^{n-1}\otimes_{\MU}\Zzp)^{hT^{n-1}}_p\to W(n-1)^{hT^{n-1}}_p$ obtained by \ref{prop:wn} and tensoring with $K(n)$, we obtain an $\Ee_0$-map 
	\[K(n-1)\otimes(T^n\otimes\Ff)^{hT^n}\to K(n-1)\otimes W(n-1)^{hT^{n-1}}_p.\]
	By Corollary \ref{cor:wnp}, the unit is nonzero in $K(n-1)_*(W(n-1)^{hT^{n-1}}_p)$. Hence the unit is also nonzero in $K(n-1)_*((T^n\otimes\Ff)^{hT^n})$. In particular, $K(n-1)\otimes(T^n\otimes\Ff)^{hT^n}$ is nonzero.
\end{proof} 

\begin{proof}[Proof of Theorem \ref{thm:main2}]
	Since the height of $\Zz$ is zero, \cite[Theorem E]{BSY22} implies that the height of $K^{(n)}(\Zz)$ is $n$. The trace map gives an $\Ee_1$-map 
	\[0\simeq K(n+1)\otimes K^{(n)}(\Zz)\longrightarrow K(n+1)\otimes(T^n\otimes\Zz)^{hT^n},\]
	which implies that $K(n+1)\otimes(T^n\otimes\Zz)^{hT^n}\simeq 0$. On the other hand, the $\Ee_{\infty}$-map \eqref{eq:forget} and Proposition \ref{prop:wn} give an $\Ee_1$-map 
	\[K(n)\otimes(T^n\otimes\Zz)^{hT^n}\longrightarrow K(n)\otimes W(n)^{hT^n}.\]
	By Proposition \ref{prop:k(n)wn}, $K(n)\otimes W(n)^{hT^n}$ is nonzero, and so is $K(n)\otimes(T^n\otimes\Zz)^{hT^n}$.
\end{proof}

Next we prove Theorem \ref{thm:main3}. For this, let $\ell$ be the Adams summand of $ku_{(p)}$. Recall that, for an odd prime $p$, $j_p$ is defined to be $\tau_{\geq0}L_{K(1)}\Ss$. Applying $L_{K(1)}$ to the unit map $\Ss\to KU$, we obtain an $\Ee_{\infty}$-map $L_{K(1)}\Ss\to(KU_p)^{h\Ff^{\times}}$, where the action of $\Ff^{\times}$ on $KU$ is given by the Adams operation. Applying $\tau_{\geq0}$ to both sides, we obtain an $\Ee_{\infty}$-map $j_p\to\ell_p$. When $p=2$, \cite[Remark 2.17]{LL26} implies that there is an equivlance $j_2\simeq (ko_2)^{h\Zz}\times_{(\tau_{\leq2}ko_2)^{h\Zz}}\tau_{\leq2}\Ss_2$. Hence we can also obtain an $\Ee_{\infty}$-map
\[j_2\simeq (ko_2)^{h\Zz}\times_{(\tau_{\leq2}ko_2)^{h\Zz}}\tau_{\leq2}\Ss_2\longrightarrow(ko_2)^{h\Zz}\longrightarrow ko_2\longrightarrow ku_2=\ell_2,\]
where the last map is the complexification map.

Moreover, let $\eta\colon\MU\to\ell$ be the unit map and $\iota\colon\Zzp[v_1]\to\MU_*$ be the canonical inclusion. Set $w_i=v_i-\iota\circ\eta_*(v_i)$ and $x^{\prime}_j=x_j-\iota\circ\eta_*(x_i)$ for $i\geq2$ and $j\neq p^k-1$. Since $\iota\circ\eta_*(v_i)$ and $\iota\circ\eta_*(x_i)$ are polynomials of $v_1$, there is an isomorphism
\[\MU_*\cong\Zzp[v_1,w_i,x^{\prime}_j\ |\ i\geq2,\ j\neq p^k-1]\]
with $\eta_*(w_i)=\eta_*(x^{\prime}_j)=0$. In particular, we obtain the element $z^{\prime}_i=\sigma^2w_i\in\pi_{2p^i}(T\otimes_{\MU}\ell)$ for $i\geq2$.

\begin{proof}[Proof of Theorem \ref{thm:main3}]
	Applying $K(n+1)\otimes(T^n\otimes-)^{hT^n}$ to the $\Ee_{\infty}$-map $j_p\to\ell_p$ and using the $\Ee_{\infty}$-map \eqref{eq:forget} give an $\Ee_1$-map
	\[K(n+1)\otimes(T^n\otimes j_p)^{hT^n}\longrightarrow K(n+1)\otimes(T^n\otimes_{\MU}\ell_p)^{hT^n}.\]
	Using the same argument as in Section \ref{sec3} with $\Zzp$ replaced by $\ell$, we obtain an even $\Ee_{\infty}$-$\MU$-algebra $W^{\prime}(n)$ with a $T^n$-action and a $T^n$-equivariant $\Ee_{\infty}$-$\MU$-algebra map $g^{\prime}_n\colon T^n\otimes_{\MU}\ell\to W^{\prime}(n)$. Moreover, $\pi_*(W^{\prime}(n))$ is a free $\ell_*$-module and $\pi_*(W^{\prime}(n)/(p,v_1))$ is a free $\Ff[\zeta^{\prime}_{1,2},\cdots,\zeta^{\prime}_{n,2}]$-module. Here, using $z^{\prime}_j$ above, $\zeta^{\prime}_{i,j}$ is defined as $\zeta_{i,j}$ with $\Zzp$ replaced by $\ell$ for $1\leq i\leq n$ and $j\geq2$. Lemma \ref{lem:gr and reg2} implies that $(p,v_1)$ is a regular sequence in $\pi_*(W^{\prime}(n)^{hT^n})$ and $\mathrm{gr}(\pi_*(W^{\prime}(n)^{hT^n}/(p,v_1)))\cong\pi_*(W^{\prime}(n)/(p,v_1))[t_1,\cdots,t_n]$. Moreover, the same argument as the proof of Proposition \ref{prop:k(n)wn} implies that
	\[v_j=\sum_{i=1}^nt_i\zeta^{\prime}_{i,j}=\sum_{i=1}^nt_i{\zeta^{\prime}_{i,2}}^{p^{j-2}}\]
	in $\mathrm{gr}(\pi_*(W^{\prime}(n)^{hT^n}/(p,v_1)))$ for $2\leq j\leq p+1$. Hence, considered as elements in $\pi_*(W^{\prime}(n)^{hT^n}/(p,v_1))$, $\{v_2,\cdots,v_{n+1}\}$ is a regular sequence by Lemma \ref{lem:reg in poly}. Consequently, the same argument as the proof of Proposition \ref{prop:k(n)wn} implies that $K(n+1)\otimes W^{\prime}(n)^{hT^n}\not\simeq0$. Using the $\Ee_1$-map $K(n+1)\otimes(T^n\otimes j_p)^{hT^n}\to K(n+1)\otimes W^{\prime}(n)^{hT^n}_p$ constructed above, we obtain that $K(n+1)\otimes(T^n\otimes j_p)^{hT^n}$ is also nonzero.
	
	It remains to show that $K(n+2)\otimes(T^n\otimes j_p)^{hT^n}$ is zero. The trace map gives an $\Ee_1$-map
	\[K(n+2)\otimes K^{(n)}(j_p)\longrightarrow K(n+2)\otimes(T^n\otimes j_p)^{hT^n}.\]
	Therefore, it suffices to show that $K(n+2)\otimes K^{(n)}(j_p)\simeq0$. When $p$ is odd, since $j_p\simeq \tau_{\geq0}L_{K(1)}\Ss,$ the height of $j_p$ is the same as the height of $L_{K(1)}\Ss$ by Lemma 2.2 and Lemma 2.3 in \cite{LMMT24}. We have that 
	\[K(m)\otimes L_{K(1)}\Ss\simeq K(m)\otimes L_{E_1}L_{K(1)}\Ss\simeq L_{E_1}K(m)\otimes L_{K(1)}\Ss\simeq0\] 
	for $m\geq2$, where $E_1$ is the height one Morava $E$-theory spectrum. Hence the height of $j_p$ is one. When $p=2$, Lemma 2.2 and Lemma 2.3 in \cite{LMMT24} implies that there is an equivalence
	\[L_{K(m)}j_2\simeq \mathrm{fib}(L_{K(m)}ko_2\xrightarrow{\psi-1}L_{K(m)}ko_2),\]
	where $\psi$ is the Adams operation. Since the height of $ko_2$ is one by the Wood cofiber sequence, the height of $j_2$ is at most one. When $m=1$, the equivalence $\mathrm{fib}(L_{K(1)}ko_2\xrightarrow{\psi-1}L_{K(1)}ko_2)\simeq L_{K(1)}\Ss$ implies that the height of $j_2$ is exactly one.
	Now, \cite[Theorem E]{BSY22} implies that $K^{(n)}(j_p)$ has height $n+1$ for all primes, which implies that $K(n+2)\otimes K^{(n)}(j_p)$ is equivalent to zero.
\end{proof}

\begin{rem}
	The same argument as the proof of Theorem \ref{thm:main2} implies that the height of $(T^n\otimes\ell)^{hT^n}$ is $n+1$, since the height of $\ell$ is one.
\end{rem}

We conclude with an example showing that the functor $(T^n\otimes -)^{hT^n}$ need not increase the height.

\begin{exmp}
	For $n,m\geq1$, let $E_m$ be the height $m$ Morava $E$-theory spectrum. Since $L_{E_m}$ is a smashing localization, $T^n\otimes E_m$ is $E_m$-local. Since the limit of $E_m$-local spectra is still $E_m$-local, $(T^n\otimes E_m)^{hT^n}$ is still $E_m$-local. Consequently, $(T^n\otimes E_m)^{hT^n}\simeq L_{E_m}(T^n\otimes E_m)^{hT^n}$ is of height $m$.
\end{exmp}

\end{document}